\documentclass[11pt]{article}

\usepackage{graphicx}
\usepackage{amsmath}
\usepackage{amsfonts}
\usepackage{epsfig}
\usepackage{color}
\usepackage{amsthm}

\newcommand{\bfm}[1]{\mbox{\boldmath $#1$}}
\newcommand{\bbeta}{\bfm{\beta}}
\newcommand{\bepsilon}{\bfm{\epsilon}}
\newcommand{\bg}{{\bf g}}

\newtheorem{theorem}{Theorem}

\newtheorem*{assumption}{Assumption}
\newtheorem{corollary}{Corollary}

\newtheorem{remark}{Remark}
\newtheorem{definition}{Definition}

\begin{document}

\title{\bf MAP MODEL SELECTION IN GAUSSIAN REGRESSION}

\author{{\bf Felix Abramovich}\\
Department of Statistics $\&$ Operations Research \\
Tel Aviv University \\
Tel Aviv 69978,  Israel \\ \\
{\bf Vadim Grinshtein}\\
Department of Mathematics \\
The Open University of Israel \\
Raanana 43107, Israel}

\date{}

\maketitle

\pagestyle{empty}

\begin{abstract}
We consider a Bayesian approach to model selection in Gaussian linear regression,
where the number of predictors might be much larger than the number of
observations. 
From a frequentist view, the proposed procedure results in the
penalized least squares estimation with a complexity penalty associated with  a
prior on the model size. We investigate the optimality properties of the 
resulting model selector. We establish the oracle inequality and specify conditions
on the prior that imply its asymptotic minimaxity within a wide range of sparse 
and dense settings for ``nearly-orthogonal'' and ``multicollinear'' designs.

\end{abstract}

\bibliographystyle{plain}

\section{Introduction} \label{sec:intr}
Consider the standard Gaussian linear regression model 
\begin{equation}
{\bf y}=X\bbeta+\bepsilon, \label{eq:model}
\end{equation}
where ${\bf y} \in \mathbb{R}^n$ is a vector of the observed response variable 
$Y$, $X_{n \times p}$ is the design matrix of the $p$ explanatory variables 
(predictors)  $X_1,...,X_p$,
$\bbeta \in \mathbb{R}^p$ is a vector of unknown regression coefficients, $\bepsilon \sim N({\bf 0},\sigma^2I_n)$ and the noise variance $\sigma^2$ is assumed to be known.  

A variety of statistical applications of regression models involves a vast number of potential explanatory  variables that
might be even large relatively to the amount of available data. It raises a severe
``curse of dimensionality'' problem. Reducing dimensionality of the model becomes therefore
crucial in the analysis of such  large data sets. 
The goal of model (or variable) selection is to select the ``best'', parsimonious subset of  predictors. 
The corresponding coefficients are then usually estimated by least squares.
The meaning of the ``best'' subset however depends on the particular aim at hand.  
One should distinguish, for example, between estimation of regression 
coefficients $\bbeta$, estimation of the mean vector $X\bbeta$, model 
identification and predicting future observations.  Different aims may lead to 
different optimal model selection procedures especially when the number of
potential predictors $p$ might be much larger than the sample size $n$. 
In this paper we focus
on estimating the mean vector $X\bbeta$ and the goodness of a model (subset
of predictors) $M$ is measured by the quadratic risk 
$E||X\hat{\bbeta}_M-X\bbeta||^2=||X\bbeta_M-X\bbeta||^2+\sigma^2|M|$,
where
$\hat{\bbeta}_M$ is the least squares estimate of $\bbeta$ and
$X\bbeta_M$ is the projection of $X\bbeta$ on the span of $M$. 
The first (bias) term of the risk decomposition represents the approximation error of the projection,
while the second (variance) term is the price for estimating the projection
coefficients $\bbeta_M$ by $\hat{\bbeta}_M$ and is proportional to the model 
size.
The ``best'' model then is the one with the minimal quadratic risk. 
Note that the true underlying model in (\ref{eq:model}) is not necessarily
the best in this sense since sometimes it is possible to reduce its
risk by excluding predictors with small (but still nonzero!) coefficients.

Such a criterion for model selection is obviously impossible to 
implement since it depends on the unknown $\bbeta$. 
Instead, the corresponding ideal minimal risk can be used as a benchmark for any
available model selection procedure. The model selection criteria are typically 
based on the {\em empirical} quadratic risk $||{\bf y}-X\hat{\bbeta}_M||^2$,
which is essentially the least squares.
However, direct minimization of the empirical risk evidently leads to
a trivial (unsatisfactory!) choice of the saturated model. 
A typical remedy is then to add a complexity penalty $Pen(|M|)$ 
that increases with the model size, and to consider {\em penalized} least squares criterion of the form
\begin{equation}
||{\bf y}-X\hat{\bbeta}_M||^2+Pen(|M|) \rightarrow \min_M  \label{eq:pmle}
\end{equation}

The properties of the resulting estimator depends  on the proper choice of the complexity penalty function $Pen(\cdot)$ in (\ref{eq:pmle}). 
There exists a plethora of works in literature on this problem. 
The standard, most commonly used choice is a {\em linear} type penalty of the form 
$Pen(k)=2\sigma^2 \lambda k$ for some fixed $\lambda>0$. The most known examples motivated by different ideas include AIC  for $\lambda=1$ (Akaike, 1974),  BIC  for $\lambda=(\ln n)/2$
(Schwarz, 1978) and RIC for $\lambda=\ln p$ (Foster \& George, 1994). 
A series of recent works suggested the so-called $2k\ln(p/k)$-type 
nonlinear penalties of the form
$Pen(k) =2\sigma^2c k(\ln(p/k)+\zeta_{p,k})$, where $c>1$ and $\zeta_{p,k}$ is 
some ``negligible'' term (see, e.g., Birg\'e \& Massart, 2001, 2007; 
Johnstone, 2002; Abramovich {\em et al.}, 2006; Bunea, Tsybakov \& Wegkamp, 2007).

In this paper we present a Bayesian formalism to the model selection problem in 
Gaussian linear regression (\ref{eq:model}) that leads to a general penalized 
model selection rule (\ref{eq:pmle}). The proposed
Bayesian approach can be used, in fact,
as a natural tool for obtaining a variety of penalized least
squares estimators with different complexity penalties that accommodate 
many of the known model selection procedures as particular 
cases corresponding to specific choices of the prior. 
Within Bayesian framework, the penalty term in (\ref{eq:pmle}) is interpreted 
as proportional to the logarithm of a prior distribution. 
Complexity penalties $Pen(|M|)$ imply placing a prior on the model size
(the number of nonzero entries of $\bbeta$).
Minimization of 
(\ref{eq:pmle}) corresponds to the maximum {\em a posteriori} (MAP)
rule yielding the resulting MAP model selector to be the posterior mode.

Although there exists a large amount of literature on Bayesian model selection
(see George \& McCulloch, 1993, 1997;
Chipman, George \& McCullogh, 2001; Liang {\em et al.}, 2008 for
surveys), it mainly focuses on ``purely Bayesian'' issues (e.g., prior 
specification, posterior calculations, etc.) and does not investigate
the optimality of the resulting 
Bayesian procedures from a frequentist view.
In this paper we study the optimality properties of the proposed MAP 
model selectors for estimating the mean vector $X\bbeta$ in (\ref{eq:model}).
First, under mild conditions on the prior 
we establish the oracle inequality and show that, up to a constant multiplier,
they achieve the minimal possible risk among all estimators.
We then investigate their asymptotic minimaxity. 
For ``nearly-orthogonal'' design they are proved to be
simultaneously rate-optimal (in the minimax sense) over a wide range of sparse 
and dense 
settings and outperform various existing model selection procedures, 
e.g. AIC, BIC, RIC, Lasso
(Tibshirani, 1996) and Dantzig selector (Cand\'es \& Tao, 2007). 
In a way, these results extend those
of Abramovich, Grinshtein \& Pensky (2007) and Abramovich {\em et al.} (2010)
for the normal means problem corresponding to the particular case $X=I_n$.

The analysis of ``multicollinear'' design, which is especially relevant
for ``$p$ much larger than $n$'' setup, is more delicate.
We demonstrate that the lower bounds for the minimax rates for
estimating the mean vector in this case are smaller than those for 
``nearly-orthogonal'' design by the factor
depending on the design properties. Such ``blessing of multicollinearity'' can
be explained by a possibility of exploiting correlations between predictors
to reduce the size of a model (hence, to decrease the variance) without paying
much extra price in the bias term.
We show that 
under some additional assumptions on the design and the 
coefficients vector $\bbeta$ in (\ref{eq:model}), the proposed Bayesian
model selectors are still asymptotically rate-optimal.

The paper is organized as follows. 
The Bayesian model selection procedure that leads to a penalized least squares
estimator (\ref{eq:pmle}) is introduced in Section \ref{sec:MAP}. 
In Section \ref{sec:oracle} we derive an
upper bound for the quadratic risk of the resulting MAP model selector, compare it
with that of an oracle and find the conditions on the prior where, up to
a constant multiplier, it achieves the minimal possible risk among all estimators. 
In Section \ref{sec:minimax} we obtain the upper and lower risk bounds of the 
MAP model selector in a sparse setup that allows us to investigate its
asymptotic minimaxity for nearly-orthogonal and multicollinear designs in
Section \ref{sec:asymp}. The computational aspects are discussed in Section
\ref{sec:comp}, and 
the main take-away messages of the paper are
summarized in Section \ref{sec:disc}. All the proofs are given in the Appendix.

\section{MAP model selection procedure} \label{sec:MAP}
Consider the Gaussian linear regression model (\ref{eq:model}), where
the number of possible predictors $p$ might be even larger then the number of
observations $n$. Let $r=rank(X) (\leq \min(p,n))$ and assume that any 
$r$ columns of $X$ are linearly independent.
For the ``standard'' linear regression setup, where all $p$
predictors are linearly independent and there are at least $p$ linearly
independent design points, $r=p$.

Any model $M$ is uniquely  defined by the $p \times p$ diagonal indicator matrix
$D_M=diag({\bf d}_M)$, where $d_{jM}=\mathbb{I}\{X_j \in M\}$ and, therefore,
$|M|=tr(D_M)$. The corresponding least square estimate
$\hat{\bbeta}_M=(D_M X'X D_M)^+ D_M X' \bf{y}$, where ``+'' denotes the
generalized inverse matrix. 

Assume some prior on the model size $\pi(k)=P(|M|=k)$ , 
where $\pi(k)>0,\;k=0,...,r$ ($k=0$ corresponds to a null model with a single
intercept) and $\pi(k)=0$ for $k>r$ since 
otherwise, there necessarily exists another vector 
$\bbeta^*$ with at most $r$ nonzero
entries also satisfying (\ref{eq:model}), that is, $X\bbeta=X\bbeta^*$.

For any $k=0,...,r-1$ there are ${p \choose k}$ different models of a given size $k$. Assume  all of them to be equally
likely, that is, conditionally on $|M|=k$,
$$
P(M\; \bigl|\; |M|=k) = {p \choose k}^{-1}
$$
One should be a little bit more careful for $k=r=rank(X)$. Although there are 
${p \choose r}$ different sets of predictors of size $r$, all of them 
evidently result in the same estimator for the mean vector and, in this sense, 
are essentially undistinguishable and associated with a {\em single} (saturated)
model. Hence, in this case, we set 
\begin{equation}
P(M\; \bigl|\; |M|=r)=1 \label{eq:full}
\end{equation}
Finally, assume the normal prior on the unknown vector of $k$ coefficients of the model $M$:
$\bbeta_M \sim N_p({\bf 0},\gamma \sigma^2 (D_M X'XD_M)^{+})$. 
This is a well-known conventional $g$-prior of Zellner (1986). 
%See Liang {\em et al.} (2008) for a survey on $g$-prior and related priors on
%$\bbeta_M$ within Bayesian perspective.
  
For the proposed hierarchical prior, straightforward calculus yields the posterior probability of a model $M$ of size $|M|=0,...,r-1$~:
\begin{equation}
P(M|\bfm{y}) \propto \pi(|M|) {p \choose |M|}^{-1} (1+\gamma)^{-\frac{|M|}{2}}
\exp\left\{\frac{\gamma}{\gamma+1}\frac{{\bf y}'XD_M(D_M X'XD_M)^{+}D_MX'{\bf y}}{2\sigma^2}\right\} \label{eq:post}
\end{equation}
Finding the most likely model leads therefore to the following maximum {\em 
a posteriori} (MAP) model selection criterion:
$$
{\bf y}'XD_M(D_M X'XD_M)^{+}D_MX'{\bf y}+2\sigma^2(1+1/\gamma)\ln\left\{{p\choose |M|}^{-1}\pi(|M|)(1+\gamma)^{-\frac{|M|}{2}}\right\} \rightarrow \max_M
$$
or, equivalently,
\begin{equation}
||{\bf y}-X\hat{\bbeta}_M||^2+2\sigma^2(1+1/\gamma)\ln\left\{{p\choose |M|}\pi(|M|)^{-1}(1+\gamma)^{\frac{|M|}{2}}\right\} \rightarrow \min_M,  \label{eq:map}
\end{equation}
which is of the general type (\ref{eq:pmle}) with the complexity penalty
\begin{equation}
Pen(k)=2\sigma^2(1+1/\gamma)\ln\left\{{p\choose k}\pi(k)^{-1}(1+\gamma)^{\frac{k}{2}}\right\},\;\;\;k=0,...,r-1 \label{eq:pen}
\end{equation}
Similarly, for $|M|=r$ from (\ref{eq:full}) one has
\begin{equation}
Pen(r)=2\sigma^2(1+1/\gamma)\ln\left\{\pi(r)^{-1}(1+\gamma)^{\frac{r}{2}}\right\} \label{eq:pen1}
\end{equation}

A specific form of the penalty (\ref{eq:pen})-(\ref{eq:pen1}) 
depends on the choice of a prior $\pi(\cdot)$.  
In particular, the (truncated if $p>n$) binomial prior $B(p,\xi)$ 
%, where $\pi(k) \propto {p \choose k} \xi^k(1-\xi)^{p-k},\;k=0,...,r$, 
corresponds to the prior assumption that the indicators $d_{jM}$ are independent. 
The binomial prior yields the linear penalty $Pen(k)=2\sigma^2\lambda k$, where $\lambda=(1+1/\gamma)\ln\{\sqrt{1+\gamma}(1-\xi)/\xi\}
\sim \ln\{\sqrt{\gamma}(1-\xi)/\xi\}$ for sufficiently large variance ratio $\gamma$. The AIC criterion corresponds then to 
$\xi \sim \sqrt{\gamma}/(e+\sqrt{\gamma})$, while 
$\xi \sim \sqrt{\gamma}/(p+\sqrt{\gamma})$ leads to the RIC criterion.  
These relations indicate that RIC should be appropriate for sparse cases,
where the size of the true (unknown) model is believed to be much less 
than the number of possible predictors, while AIC is suitable for dense cases, 
where they are of the same order. In fact, any binomial prior or, equivalently, any linear penalty cannot ``kill two birds with one stone''. 
On the other hand, the (truncated) geometric prior $\pi(k) \propto q^k,\;k=1,...,r$ for some $0 < q < 1$, implies 
$Pen(k) \sim 2\sigma^2(1+1/\gamma)k(\ln(p/k)+c(\gamma,q))$ which is of the $2k\ln(p/k)$-type introduced above. For large $\gamma$ it behaves similar to RIC for $k \ll p$ and to AIC for $k \sim p$
and is, therefore, adaptive to both sparse and dense cases. We will discuss these issues more rigorously in Section \ref{sec:asymp} below.
 
\section{Oracle inequality} \label{sec:oracle}
In this section we derive an upper bound for the quadratic risk of the proposed MAP model selector and compare it with the ideal minimal quadratic risk often 
called in literature as an oracle risk.

\begin{assumption}[P] \label{as:P}
Assume that 
$$
\pi(k) \leq {p \choose k}e^{-c(\gamma)k},\;k=0,...,r-1,\;{\rm and} \;\pi(r) \leq e^{-c(\gamma)r}
$$
where $c(\gamma)=8(\gamma+3/4)^2 > 9/2$. 
\end{assumption}
Assumption (P) is not restrictive. Indeed, the obvious inequality 
${p \choose k} \geq (p/k)^k$ implies that for $k < r$ it automatically holds for{\em any} prior
$\pi(k)$ for all $k \leq pe^{-c(\gamma)}$. Assumption (P) is used to 
establish an upper bound for the quadratic risk of the MAP model selector.

\begin{theorem} \label{th:bound}
Let the model $\hat{M}$ be the solution of (\ref{eq:pmle}) with the
complexity penalty $Pen(\cdot)$ given in (\ref{eq:pen})-(\ref{eq:pen1}) and 
$\hat{\bbeta}_{\hat M}$ be the corresponding least squares estimate.
Then, under Assumption (P)
\begin{equation}
E||X\hat{\bbeta}_{\hat M}-X\bbeta||^2 \leq c_0(\gamma)
\inf_M\left\{||X\bbeta_M-X\bbeta||^2+Pen(|M|)\right\}+c_1(\gamma)\sigma^2 
\label{eq:bound}
\end{equation}
for some $c_0(\gamma)$ and $c_1(\gamma)$ depending only on $\gamma$. 

%\begin{eqnarray}
%E||X\hat{\bbeta}_{\hat M}-X\bbeta||^2 & \leq & c_0(\gamma)
%\inf_M\left\{||X\bbeta_M-X\bbeta||^2+
%2\sigma^2(1+1/\gamma) \ln\left({p \choose |M|}\pi(|M|)^{-1}(1+\gamma)^{\frac{|M|}{2}}\right)\right\} \nonumber \\ 
%& + & c_1(\gamma)\sigma^2 \nonumber \\ \label{eq:bound}
%\end{eqnarray}
%for some $c_0(\gamma)$ and $c_1(\gamma)$ depending only on $\gamma$.

\end{theorem}

To assess the quality of the upper bound in (\ref{eq:bound}), we compare
it with the oracle risk $\inf_M E||X\hat{\bbeta}_M-X\bbeta||^2$. 
Note that
the oracle risk is exactly zero when $\bbeta \equiv {\bf 0}$ and, evidently,
no estimator can achieve it in this case. Hence, an additional, typically 
negligible
term $\sigma^2$, which is, essentially, an error of estimating a single
extra parameter, is usually added to the oracle risk for a proper comparison.
It is known that no estimator can attain a risk smaller than
within $2\ln p$ factor from that of an oracle (e.g., Foster \& George, 1994;
Donoho \& Johnstone, 1995; Cand\`es, 2006). The following theorem shows that under certain
additional conditions on the prior $\pi(\cdot)$, the resulting MAP model 
selector achieves this minimal possible risk among all estimators up to 
a constant multiplier depending on $\gamma$:
\begin{theorem}[oracle inequality] \label{th:oracle}
Let $\pi(k)$ satisfy Assumption (P) and, in addition, 
$\pi(0) \geq p^{-c},\; 
\pi(k) \geq  p^{-ck},\;k=1,...,r$ for some constant $c>0$. 
Then, the resulting MAP model selector satisfies
$$
E||X\hat{\bbeta}_{\hat M}-X\bbeta||^2 \leq c_2(\gamma)\ln p 
\left(\inf_M E||X\hat{\bbeta}_M-X\bbeta||^2 + \sigma^2\right)
$$
for some $c_2(\gamma) \geq 1$.
\end{theorem}
In particular, it can be easily shown that Theorem \ref{th:oracle} holds for the
(truncated) binomial $B(p,\xi)$ with $\xi=1/p$ (RIC criterion)
and geometric priors (see Section \ref{sec:MAP}). More generally,
all priors such that
$\ln \pi(k)=O(k\ln(k/p))$ corresponding to the 
$2k\ln(p/k)$-type penalties satisfy the conditions of Theorem \ref{th:oracle}. 

\section{Risk bounds for sparse settings} \label{sec:minimax}
In the previous section we considered the global behavior of the MAP
estimator without any restrictions on the model size. However, in the
analysis of large data sets, it is typically reasonable to assume that 
the true model in (\ref{eq:model}) is sparse in the sense that only part of 
coefficients in $\bbeta$ are different from zero.
We now show that under such extra sparsity assumption, more can be said on the 
optimality of the MAP model selection.
  
For a given $1 \leq p_0 \leq r$, define the sets
of models ${\cal M}_{p_0}$ that have at most $p_0$ predictors, that is,
${\cal M}_{p_0}=\{M : |M| \leq p_0\}$. Obviously, if a true model in
(\ref{eq:model}) belongs to ${\cal M}_{p_0}$, the $l_0$ quasi-norm of the
corresponding coefficients vector $||\bbeta||_0 \leq p_0$, where $||\bbeta||_0$
is the number of its nonzero entries.
In this section we find the upper and lower bounds for the maximal risk of 
the proposed MAP model selector over ${\cal M}_{p_0}$. 

\begin{theorem} \label{th:upper} 
Let $1 \leq p_0  \leq r$, 
$\pi(\cdot)$ satisfy Assumption (P) and, in addition, 
$\pi(p_0) \geq (p_0/(pe))^{cp_0}$ if $p_0 < r$ or $\pi(r) \geq e^{-cr}$
if $p_0=r$ for some constant $c>c(\gamma)$.
Then, there exists a constant $C_1(\gamma)>0$ depending
only on $\gamma$ such that
\begin{equation}
\sup_{\bbeta: ||\bbeta||_0 \leq p_0}E||X\hat{\bbeta}_{\hat M}-X\bbeta||^2 
\leq C_1(\gamma) \sigma^2 \min\left(p_0 (\ln(p/p_0)+1),r\right) \label{eq:upper}
\end{equation}
\end{theorem}
The general upper bound (\ref{eq:upper}) for the maximal risk of the MAP 
selector over ${\cal M}_{p_0}$
in Theorem \ref{th:upper} holds for any design matrix $X$. 
To assess its accuracy we establish the lower bound for the minimax
risk of estimating the mean vector $X\bbeta$ in (\ref{eq:model}).

For any given $k=1,...,r$, let $\phi_{min}[k]$ and $\phi_{max}[k]$
be the $k$-sparse minimal and maximal eigenvalues of the design defined as 
$$
\phi_{min}[k]=\min_{\bbeta: 1 \leq ||\bbeta||_0 \leq k} 
\frac{||X\bbeta||^2}{||\bbeta||^2},
$$
$$
\phi_{max}[k]=\max_{\bbeta: 1 \leq ||\bbeta||_0 \leq k} 
\frac{||X\bbeta||^2}{||\bbeta||^2}
$$
(see Meinshausen \& Yu, 2009; Bickel, Ritov \& Tsybakov, 2009).
In fact, $\phi_{min}[k]$ and $\phi_{max}[k]$ are respectively
the minimal and maximal eigenvalues of all $k \times k$ submatrices  
of the matrix $X'X$ generated by any $k$ columns of $X$. 
%Obviously, if $p>n$, then $\phi_{min}[k]=0$ for any $k>n$.
Let $\tau[k]=\phi_{min}[k]/\phi_{max}[k],\;k=1,...,r$ and
set $\tau[k]=\tau[r]$ for all $k > r$. 
By the definition, $\tau[k]$ is a non-increasing function of $k$.
Obviously, $\tau[k] \leq 1$ and for the orthogonal design
the equality holds for all $k$. 

\begin{theorem} \label{th:lower}
Consider the model (\ref{eq:model}) and let $1 \leq p_0 \leq r$.
There exists a universal constant $C_2>0$ such that
\begin{equation}
\label{eq:lower}
\inf_{\hat{\bf y}} \sup_{\bbeta: ||\bbeta||_0 \leq p_0} 
E||\hat{\bf y}-X\bbeta||^2 \geq  
\left\{
\begin{array}{ll}
C_2 \sigma^2 
\tau[2p_0]\; p_0 (\ln(p/p_0)+1) &,\; 1 \leq p_0 \leq r/2  \\
C_2 \sigma^2 \tau[p_0]\; r &,\; r/2 \leq p_0 \leq r 
\end{array}
\right.
\end{equation}
where the infimum is taken over all estimates $\hat{\bf y}$ of the mean vector 
$X\bbeta$.
\end{theorem}

Theorem \ref{th:lower} shows that the minimax lower bound (\ref{eq:lower}) depends
on a specific design matrix $X$ only through the sparse eigenvalues ratios.
A computationally simpler but less accurate minimax lower bound can be
obtained  by replacing $\tau[2p_0]$ and $\tau[p_0]$ in (\ref{eq:lower}) by 
$\tau[r]$, that for the case $r=p \leq  n$ is just the ratio of 
the minimal and maximal eigenvalues of $X'X$.

For the orthogonal design, where $\tau[\cdot] \equiv 1$, and
$p \leq n$ analogous results were obtained in Birg\'e \& Massart (2001).
For a general design and $p_0 \leq r/2$
similar minimax lower bounds were independently obtained in
Raskutti, Wainwright \& Yu (2009) for a design matrix of a full rank  
and in Rigollet \& Tsybakov (2010) for a general case within a related 
aggregation context.

The established upper and lower bounds (\ref{eq:upper}), (\ref{eq:lower})
for the risk of the MAP model selector allow us in the following section
to investigate its asymptotic minimaxity as both $n$ and $p$ increase.

\section{Asymptotic adaptive minimaxity} \label{sec:asymp}
\subsection{Nearly-orthogonal design} \label{subsec:orthog}
In this section we consider
the asymptotic properties of the MAP model selector as the sample size
$n$ increases. We allow $p=p_n$ to increase with $n$ as well and look for
a projection of the unknown mean vector on an expanding span of predictors. 
In particular, the
most challenging cases intensively studied nowadays in literature are those,
where $p > n$ or even $p \gg n$. 
In such asymptotic settings one should essentially consider a {\em sequence}
of design matrices $X_{n,p_n},$ where $r_n \rightarrow \infty$.
For simplicity of exposition, in what follows
we omit the index $n$ and denote $X_{n,p_n}$ by $X_p$ emphasizing the
dependence on the number of predictors $p$ and let $r$ tend to infinity.
Similarly, we consider now 
sequences of coefficients vectors $\bbeta_p$ and priors $\pi_p(\cdot)$.
In these notations, the original model (\ref{eq:model}) is transformed into a 
sequence of models
\begin{equation}
{\bf y}=X_p \bbeta_p + \bepsilon, \label{eq:model1}
\end{equation}
where $rank(X_p)=r$ and any $r$ columns of $X_p$ are linearly
independent (hence, $\tau_p[r]>0$),
$\bepsilon \sim N({\bf 0},\sigma^2 I_n)$ and
the noise variance $\sigma^2$ does not depend on $n$ and $p$.
One can also view a sequence of models (\ref{eq:model1}) in a 
triangular array setup (Greenshtein \& Ritov, 2004).

\begin{definition} \label{def:nearly} Consider the sequence of design
matrices $X_p$. 
The design is called nearly-orthogonal if the corresponding sequence
of sparse eigenvalues ratios $\tau_p[r]$
is bounded away from zero by some constant $c>0$. Otherwise, the design is called
multicollinear.
\end{definition}

Nearly-orthogonality condition essentially means that there is no 
multicollinearity in the design in the sense that 
there are no ``too strong'' linear relationships within any
set of $r$ columns of $X_p$. Intuitively, it is clear that in this case
$p$ cannot be ``too large'' relative to $r$ and, therefore, to $n$. Indeed, 
apply the upper and lower bounds (\ref{eq:upper}), (\ref{eq:lower}) for
$p_0=r/2$ to get $(C_2/2) \sigma^2 \tau[r] r (\ln(2p/r)+1) \leq C_1(\gamma)
\sigma^2 r$ that implies the following remark:
\begin{remark} \label{rem:orth}
For nearly-orthogonal design, necessarily $p=O(r)$ and,
therefore, $p=O(n)$.
\end{remark}

%In particular, for the orthogonal
%design evidently $\tau(k)=1$ for all $k$.

%In what follows $g_1(p) \asymp g_2(p)$ denotes 
%$0 < \lim \inf\{g_1(p)/g_2(p)\} \leq \lim \sup\{g_1(p)/g_2(p)\}$ as
%$r \rightarrow \infty$. 

The following corollary is an immediate consequence of Theorems \ref{th:upper} and \ref{th:lower}: 
\begin{corollary} \label{cor:minimax} Let the design be nearly-orthogonal.
\begin{enumerate} 
\item As $r$ increases, the asymptotic minimax risk of estimating the mean vector $X_p\bbeta_p$
over ${\cal M}_{p_0}$ is of the order $\min(p_0 (\ln(p/p_0)+1),r)$, that is,
there exist two constants $0< C_1 \leq C_2 < \infty$ such that for 
all sufficiently large $r$,
\begin{eqnarray*}
C_1\; \sigma^2 \min\left(p_0 (\ln(p/p_0)+1),r\right) &  \leq &
\inf_{\hat{y}}\sup_{\bbeta_p: ||\bbeta_p||_0 \leq p_0}  E||\hat{\bf y}-X_p\bbeta_p||^2 \\
& \leq &  C_2\;\sigma^2 \min\left(p_0 (\ln(p/p_0)+1),r\right)
\end{eqnarray*}
for all $1 \leq p_0 \leq r$.

\item Assume Assumption (P) and, in addition, that 
$\pi_p(k) \geq (k/(pe))^{c_1k},\;k=1,...,r-1$ and $\pi_p(r) \geq e^{-c_2r}$
for some constants $c_1, c_2> c(\gamma)$. 
Then, the corresponding MAP model selector attains the minimax convergence rates
simultaneously over all ${\cal M}_{p_0},\;1 \leq p_0 \leq r$.
\end{enumerate}
\end{corollary}

One can easily verify that the conditions on the prior of Corollary \ref{cor:minimax}
are satisfied, for example, for the truncated geometric prior 
(see Section \ref{sec:MAP}) for all $k=1,...,r$. 
The resulting MAP model selector attains, therefore, the 
minimax rates simultaneously for all ${\cal M}_{p_0},\;p_0=1,...,r$.
As we have mentioned, the corresponding penalty in (\ref{eq:pen}) is of
the $2k\ln(p/k)$-type.
On the other hand, no truncated binomial prior $B(p,\xi_p)$ 
can satisfy these conditions on the entire range $k=1,...,r$. 
It is easy to verify that they hold for small
$\xi_p$ if $k \ll r$ but for large $\xi_p$ if $k \sim r$.
In fact, these arguments go along the lines with the similar results of 
Foster \& George (1994) and Birg\'e \& Massart (2001, 2007). 
Recall that binomial prior corresponds to linear penalties
of the type $Pen(k)=2\sigma^2 \lambda k$ in (\ref{eq:pmle}) (see Section 
\ref{sec:MAP}). 
Foster \& George (1994) and Birg\'e \& Massart (2001, Section 5.2) showed that
the best possible risk of such estimators over ${\cal M}_{p_0}$ is 
only of order $\sigma^2 p_0\ln p$ achieved for $\lambda \sim \ln p$
corresponding to the RIC criterion. It is of the same
order as the optimal risk $\sigma^2 p_0(\ln(p/p_0) +1)$ for $p_0 \ll p$
(sparse case) but larger for dense case ($p_0 \sim p$). 
On the other hand, the risk of the
AIC estimator ($\lambda=1$) is of the order $\sigma^2 r$, which is optimal
for dense but much larger for sparse case.

Furthermore, under somewhat similar
nearly-orthogonality conditions, Bickel, Ritov \& Tsybakov (2009) showed that the 
well-known Lasso (Tibshirani, 1996) and Dantzig (Cand\'es \& Tao, 2007) 
estimators achieve only the same sub-optimal rate 
$\sigma^2p_0\ln p$ as RIC. These results are, in fact, not so surprising 
since both Lasso and Dantzig estimators are essentially based on convex 
relaxations of the $l_0$-norm of regression coefficients $||\bbeta||_0$ in the 
linear complexity penalty $2\sigma^2\lambda ||\bbeta||_0$ in order to
replace the original combinatorial problem (\ref{eq:pmle}) by a convex program. 
Thus, Lasso approximates the $l_0$-norm $||\bbeta||_0$ by the
the corresponding $l_1$-norm $||\bbeta||_1$. 
In particular, for the orthogonal design, 
linear complexity penalties and Lasso yield respectively hard and soft 
thresholding of components of $\bbeta$ with a {\em fixed} threshold. 
RIC estimator and Lasso with the optimally chosen  
tuning parameter (e.g., Bickel, Ritov \& Tsybakov, 2009)
result in this case in the well-known hard and soft universal 
thresholding of Donoho \& Johnstone (1994) with a fixed threshold
$\sigma\sqrt{2\ln p}$ which is rate-optimal for various 
sparse but not dense settings. 
On the other hand, the nonlinear MAP penalty corresponds to hard thresholding with
a {\em data-driven} threshold that under conditions on $\pi_p(\cdot)$ in
Corollary \ref{cor:minimax} is simultaneously minimax for both sparse and dense
cases (Abramovich, Grinshstein \& Pensky, 2007; Abramovich {\em et al.}, 2010).

Finally, note that for the nearly-orthogonal design, 
$||X_p\hat{\bbeta}_{p\hat M}-X_p\bbeta_p|| \asymp ||\hat{\bbeta}_{p\hat M}-\bbeta_p||$,
where ``$\asymp$'' means that their ratio is bounded from below and above.
Therefore, all the results of Corollary \ref{cor:minimax} for estimating the 
mean vector $X_p\bbeta_p$ in (\ref{eq:model1}) can be straightforwardly applied
for estimating the regression coefficients $\bbeta_p$. This equivalence, however, 
does not hold for the multicollinear design considered below.

\subsection{Multicollinear design} \label{subsec:multi}
Nearly-orthogonality assumption may be reasonable in the ``classical'' setup,
where $p$ is not too large relatively to $n$ but might be questionable for the analysis of
high-dimensional data, where $p \gg n$, due to the
multicollinearity phenomenon (see also Remark \ref{rem:orth}). 
When this assumption does not hold, 
the sparse eigenvalues ratios in (\ref{eq:lower}) may tend to zero as $p$ 
increases and, thus, decrease the minimax lower bound rate relatively to
the nearly-orthogonal design.
% (\ref{eq:lower}) by the factor 
%$\tau_p[2p_0]$ relatively to the nearly-orthogonal design.
In this case there is a gap between the rates in the lower and upper bounds 
(\ref{eq:lower}) and (\ref{eq:upper}). Intuitively, one can think of exploiting
correlations between predictors to reduce the size of a 
model (hence, to decrease the variance) 
without paying much extra price in the bias term, and, therefore, to reduce the risk. 
We show that under certain additional assumptions on the design and the
coefficients vector in (\ref{eq:model1}), the upper risk bound (\ref{eq:upper}) 
can be indeed reduced to the minimax lower bound  rate (\ref{eq:lower}).

For simplicity of exposition we consider the sparse case $p_0 \leq r/2$ 
although the corresponding
conditions for the dense case $r/2 \leq p_0 \leq r$ can be obtained 
in a similar way with necessary changes. 

We introduce now several definitions that will be used in the sequel
(including the proofs in the Appendix). 
For a given index set ${\cal J}$ and $l \geq |{\cal J}|$ define a 
$l \times |{\cal J}|$ matrix
$G_{l,{\cal J}}$ which columns ${\bf e}_j,\;j \in {\cal J}$ are the elements
of the standard basis in $\mathbb{R}^l$. 
Thus, for any matrix $A$ with $l$ columns,
$A G_{l,{\cal J}}$ selects the columns of $A$ indexed by ${\cal J}$.
Similarly, for any $l \times l$ symmetric matrix $A$,
$G'_{l,{\cal J}}AG_{l,{\cal J}}$ generates a (symmetric) 
$|{\cal J}| \times |{\cal J}|$ submatrix of $A$ of the corresponding
columns and rows.

For all $k=1,...,r/2$, define $k'=\lceil \tau_p[2k]\cdot k \rceil \geq 1$. 
Let ${\cal J}_M$ be an index set of predictors included in a model $M$ of size
$k=|{\cal J}_M|$. 
For any submodel $M' \subset M$ of
size $k'<k$ let $\tilde{\phi}_{M,M'}$ be the minimal eigenvalue of the
$(k-k') \times (k-k')$ matrix $\Lambda_p(M,M')=G'_{k,{\cal J}_{M/M'}}(G_{p,{\cal J}_M}'X_p'X_pG_{p,{\cal J}_M})^{-1}G_{k,{\cal J}_{M/M'}}$.
In fact, $\sigma^2\Lambda_p(M,M')$ is the covariance matrix of the components of
the least squares estimate vector $\hat{\bbeta}_M$ corresponding to 
a subset of predictors in $M/M'$.

Finally, define 
$$
\tilde{\phi}_p[k]=\min_{M: |M|=k} \max_{M' \subset M:  |M'|=k'} 
\tilde{\phi}_{M,M'}
%\tilde{\phi}_{min}(G'_{M/M'}(G_M'X_p'X_pG_M)^{-1}G_{M/M'})
$$
As we show later (see the proof of Theorem \ref{th:nonorth} in the Appendix),
$\tilde{\phi}^{-1}_p[k]$ measures an error of approximating mean vectors 
$X_p\bbeta_p$, where 
$||\bbeta_p||_0=k$, by their projections on lower dimensional subspans of
predictors. The stronger is multicollinearity, the better is the approximation
and the larger is $\tilde{\phi}_p[k]$. 

\begin{theorem} \label{th:nonorth} 
Let $\tau_p[r] \rightarrow 0$ as $r \rightarrow \infty$
(multicollinear design).
Assume the following additional assumptions on the design matrix $X_p$ and the 
(unknown) vector of
coefficients $\bbeta_p$ in (\ref{eq:model1}):
\begin{description}
\item[(D)] for all $p$ there exist $1 \leq \kappa_{p1} \leq \kappa_{p2} \leq r/2$ such that
\begin{enumerate} 
\item $\tilde{c}_1 \leq \tau_p[2k]\cdot k \leq k-1,\;k=\kappa_{p1},...,
\kappa_{p2}$
\item $\tau_p[2\kappa_{p2}] \geq (\kappa_{p2}/(pe))^{\tilde{c}_2}$
\item $\phi_{p,min}[2k]\cdot \tilde{\phi}_p[k] \geq 
\tilde{c}_3,\;k=\kappa_{p1},...,\kappa_{p2}$
\end{enumerate}

\item[(B)] $||\bbeta_p||^2_\infty \leq \tilde{c}_4 
\tau_p[2p_0] \cdot \tilde{\phi}_p[p_0]\cdot(\ln(p/p_0)+1)$, where
$p_0=||\bbeta_p||_0$
\end{description}
for some positive constants $\tilde{c}_1,\;\tilde{c}_2,\;\tilde{c}_3$
and $\tilde{c}_4$.

Then, under the above additional restrictions, if
the prior $\pi_p(\cdot)$ satisfies Assumption (P)
and for all $k=\kappa_{p1},...,\kappa_{p2}$,
$\pi_p(k') \geq (k'/(pe))^{ck'}$ for some 
positive $c>c(\gamma)$, where $k'=\lceil \tau_p[2k]\cdot k \rceil$,
the corresponding MAP model 
selector is asymptotically simultaneously  minimax 
(up to a constant multiplier) over all 
${\cal M}_{p_0},\;\;\kappa_{p1} \leq p_0 \leq \kappa_{p2}$.
\end{theorem}
Note that by simple algebra one can verify that 
$\phi_{p,min}[2k]\cdot \tilde{\phi}_p[k]
\leq 1$ and, therefore, the constant $\tilde{c}_3$ in Assumption (D.3) is not
larger than one. 

We have argued that multicollinearity typically arises when $p \gg n$.
One can easily verify that for $n=O(p^{\alpha}),\;0 \leq \alpha<1$, Assumption (D.2)
always follows from Assumption (D.1) and, therefore, can be omitted in this
case.

As we show in the proof, Assumptions (D.1, D.2) and
Assumption (B) allow one to reduce the upper bound
(\ref{eq:upper}) for the risk of the MAP model selector by the factor 
$\tau_p[2p_0]$,
while Assumption (D.3) is required to guarantee that the additional 
constraint on $\bbeta_p$ in Assumption (B) does not affect the lower bound 
(\ref{eq:lower}). 

To obtain asymptotic minimaxity of the MAP selector within the entire
range $1 \leq p_0 \leq r/2$ similar to Corollary \ref{cor:minimax} for the
nearly-orthogonal case, Assumptions (D) on the design
matrix are required to be satisfied for all $k=1,...,r/2$ that might be quite restrictive. 
However, the results of Theorem \ref{th:nonorth} are more general and 
show the tradeoff between relaxation of Assumptions (D) to a smaller range
of $k$ and the corresponding constriction of the adaptivity range for $p_0$.

%To give some insight about Assumptions (D) on the design, recall that
%$\tau_p[k]$ is a non-increasing function of $k$ and 
%multicollinearity assumes that $\tau_p[r] \rightarrow 0$.
%For those (small) $k$ in Theorem \ref{th:nonorth}, 
%where $\tau_p[2k]$ still does not tend to zero,
%multicollinearity does not play an important role yet: Assumptions (D) are
%quite obvious and the minimax rates over 
%corresponding ${\cal M}_k$ are of the same order as for the nearly-orthogonal
%case. The most interesting $k$ are those, where $\tau_p[2k] \rightarrow 0$ and
%the effect of multicollinearity appears. Assumption (D.1) ensures that
%$\tau_p[2k]$ cannot tend to zero faster than $k^{-1}$, while
%Assumption (D.3) implies in this case 
%$\tilde{\phi}_p[k] \rightarrow \infty$ or, 
%equivalently, $\tilde{\phi}^{-1}_p[k] \rightarrow 0$.
%From the proof of Theorem \ref{th:nonorth} it follows that
%$\tilde{\phi}^{-1}_p[k]$ measures an error of approximating mean vectors 
%$X_p\bbeta_p$, where 
%$||\bbeta_p||_0=k$, by their projections on lower dimensional subspans of
%predictors. 
%This approximation error then uniformly decreases as 
%$\tilde{\phi}^{-1}_p[k] \rightarrow 0$. It essentially
%means that the multicollinearity of design is ``uniformly spread''  among all
%predictors rather than ``localized'' within their certain subset.   

\section{Computational aspects} \label{sec:comp}
In practice, minimizing (\ref{eq:pmle}) (and (\ref{eq:map}) in particular) 
requires generally an NP-hard combinatorial search
over all possible models. During the last decade there have been substantial
efforts to develop various {\em approximated} algorithms for solving (\ref{eq:pmle})
that are computationally feasible for high-dimensional data (see, e.g.
Tropp \& Wright, 2010 for a survey and references therein). 
The common remedies involve either
greedy algorithms (e.g., forward selection, matching pursuit) approximating the
global solution by a stepwise sequence of local ones, or 
convex relaxation methods replacing the original combinatorial problem by a 
related convex program (e.g., Lasso and Dantzig selector for linear penalties).
The proposed Bayesian formalism allows one instead to use a stochastic search
variable selection (SSVS) techniques originated in 
George \& McCulloch (1993, 1997) for solving (\ref{eq:map})
by generating a sequence of models from
the posterior distribution $P(M|{\bf y})$ in (\ref{eq:post}). 
The key point is that the relevant
models with the highest posterior probabilities will appear most frequently
and can be identified even for a generated sample of a relatively small 
size avoiding computations of the entire posterior distribution.

The SSVS algorithm for the problem at hand can be basically described as follows.
As we have mentioned in Section \ref{sec:MAP}, every model $M$ is uniquely
defined by the corresponding indicator vector ${\bf d}_M$ and the joint 
posterior distribution of ${\bf d}_M$ is given by (\ref{eq:post}) (up to a 
normalizing constant). SSVS uses the Gibbs sampler to generate a sequence
of indicator vectors ${\bf d}_1,...,{\bf d}_m$ {\em componentwise} by sampling
consecutively from the conditional distributions $d_j|({\bf d}_{(-j)},{\bf y})\;,
j=1,...,p$, where ${\bf d}_{(-j)}=(d_1,...,d_{j-1},d_{j+1},...,d_p)'$.
The components $d_j$ can be trivially obtained as simulations of Bernoulli
draws,
%Note that $d_j|({\bf d}_{(-j)},{\bf y})$ are Bernoulli random variables, 
% \sim B(1,p(d_j=1|{\bf d}_{(-j)},{\bf y})$,
where from (\ref{eq:post}) the corresponding posterior odds ratio 
\begin{eqnarray*} 
\frac{P(d_j=1|{\bf d}_{(-j)},{\bf y})}{P(d_j=0|{\bf d}_{(-j)},{\bf y})} &  = &
\frac{P(d_j=1,{\bf d}_{(-j)}|{\bf y})}{P(d_j=0,{\bf d}_{(-j)}|{\bf y})} \\
& = & \frac{\pi(|d_{(-j)}|+1)}{\pi(|d_{(-j)}|)}\; \frac{|d_{(-j)}|+1}{p-|d_{(-j)}|}
\;(1+\gamma)^{-\frac{1}{2}}\; \exp\left\{\frac{\gamma}{\gamma+1}
\frac{\Delta RSS_j}{2\sigma^2}\right\}
\end{eqnarray*}
and $\Delta RSS_j>0$ is the increment in the residual sum of squares (RSS)
after dropping the $j$-th predictor from the model 
$(d_j=1,{\bf d}_{(-j)})'$. 
The resulting Gibbs sampler is computationally efficient and,
as $m$ increases, 
the empirical distribution of the generated 
sample converges to the actual posterior distribution of 
${\bf d}|{\bf y}$. After the sequence has reached
approximate stationarity, one can identify the most frequently appeared
vector(s) ${\bf d}$ as potential candidate(s) to solve (\ref{eq:map}).

%In particular, for the (truncated) binomial $B(p,\xi)$ and geometric 
%$Geom(1-q)$ priors, the corresponding posterior odds ratio become
%$$
%\frac{P(d_j=1|{\bf d}_{(-j)},{\bf y})}{P(d_j=0|{\bf d}_{(-j)},{\bf y})} =
%\frac{\xi}{1-\xi} \;(1+\gamma)^{-\frac{1}{2}}\; \exp\left\{\frac{\gamma}{\gamma+1}
%\frac{\Delta RSS_j}{2\sigma^2}\right\}
%$$
%and
%$$
%\frac{P(d_j=1|{\bf d}_{(-j)},{\bf y})}{P(d_j=0|{\bf d}_{(-j)},{\bf y})} =
%q \; \frac{|d_{(-j)}|+1}{p-|d_{(-j)}|} \;(1+\gamma)^{-\frac{1}{2}}\; \exp\left\{\frac{\gamma}{\gamma+1}
%\frac{\Delta RSS_j}{2\sigma^2}\right\}
%$$
%respectively.

\section{Concluding remarks} \label{sec:disc}
In this paper we considered a Bayesian approach to model selection in Gaussian
linear regression. From a frequentist view, the resulting MAP model selector is
a penalized least squares estimator with a complexity penalty associated with
a prior $\pi(\cdot)$ on the model size.
Although the proposed estimator was originated within Bayesian
framework, the latter was used as a natural tool to obtain a wide class of 
penalized least squares estimators with various complexity penalties. 
Thus, we believe that the main take-away messages of the paper
summarized below are of a more general interest. 
 
The first main take-away message is that neither linear complexity penalties 
(e.g., AIC, BIC and RIC) corresponding to binomial priors $\pi(\cdot)$,
nor closely related Lasso and Dantzig estimators can be
simultaneously minimax for both sparse and dense cases. 
We specify the class of priors and associated nonlinear penalties  
that do yield such a wide adaptivity range. In particular,
it includes $2k\ln(p/k)$-type penalties.

Another important take-away message is about the effect of 
multicollinearity of design.
Unlike model identification or coefficients estimation, where multicollinearity
is a ``curse'', it may become a ``blessing'' for estimating the mean vector
allowing one to exploit correlations between predictors
to reduce the size of a model (hence, to decrease the variance) without paying
much extra price in the bias term. Interestingly, a similar phenomenon occurs
in a testing setup (e.g., Hall \& Jin, 2010).

\section*{Acknowledgments}
The authors would like to thank Alexander Samarov and Ya\'acov Ritov for 
valuable remarks. The authors are especially grateful to an anonymous 
referee for an excellent constructive review of the first version of the paper.

\section{Appendix}
Throughout the proofs we use $C$ to denote a generic positive constant, not
necessarily the same each time it is used, even within a single equation.

\subsection{Proof of Theorem \ref{th:bound}}
Define 
\begin{equation}
L_k=(1/k)\ln\left({p \choose k}\pi^{-1}(k)\right) \geq c(\gamma),
\;k=1,...,r-1 \label{eq:Lk}
\end{equation}
and 
\begin{equation}
L_r=(1/r)\ln\pi^{-1}(r) \geq c(\gamma) \label{eq:Lr}
\end{equation}
In terms of $L_k$ the complexity penalty (\ref{eq:pen})-(\ref{eq:pen1}) is
$Pen(k)=\sigma^2(1+1/\gamma)k(2L_k+\ln(1+\gamma))$.
Following the arguments of the proof of Theorem 1 of Abramovich {\em et al.}
(2007), under the Assumption (P) one has
$$
\sum_{k=1}^{r-1}{p \choose k}e^{-kL_k}+e^{-rL_r} = \sum_{k=1}^r \pi(k) =
1-\pi(0) < 1
$$
and
$$
(1+1/\gamma)(2L_k+\ln(1+\gamma)) \geq C (1+\sqrt{2L_k})^2,\;k=1,...,r
$$
The proof of Theorem \ref{th:bound} then follows directly from 
Theorem 2 of Birg\'e \& Massart (2001).
\newline $\Box$

\subsection{Proof of Theorem \ref{th:oracle}}
Let $L^*=\max_{0 \leq k \leq r} L_k$, where $L_k,\;k=1,...,r$
were defined in (\ref{eq:Lk})-(\ref{eq:Lr}) and $L_0=2\ln\pi^{-1}(0)$. 
Simple calculus shows that the conditions on 
$\pi(\cdot)$ in Theorem \ref{th:oracle} imply $L^*=O(\ln p)$.
 
Consider first the case $k \geq 1$.
From Theorem \ref{th:bound} we have
\begin{eqnarray}
E||X\hat{\bbeta}_{\hat M}-X\bbeta||^2 & \leq & c_0(\gamma)
\inf_M\left\{||X\bbeta_M-X\bbeta||^2+
\sigma^2(1+1/\gamma)|M|(2L_{|M|}+\ln(1+\gamma))\right\}+c_1(\gamma)\sigma^2 \nonumber \\
& \leq & c_0(\gamma)(1+1/\gamma)(2L^*+\ln(1+\gamma))
\inf_M\left\{||X\bbeta_M-X\bbeta||^2+|M|\sigma^2\right\}+c_1(\gamma)\sigma^2
\nonumber \\ 
& \leq & c_2(\gamma)(2L^*+\ln(1+\gamma))\left\{\inf_M E||X\hat{\bbeta}_M-X\bbeta||^2 + \sigma^2\right\} \nonumber \\ \label{eq:Lk1}  
\end{eqnarray}

For the degenerative case $k=0$ ($M=\{0\}$), Theorem \ref{th:bound} implies
\begin{eqnarray}
E||X\hat{\bbeta}_{\hat M}-X\bbeta||^2 & \leq & 
%c_0(\gamma)\left\(||X\bbeta||^2+2\sigma^2(1+1/\gamma)\ln\pi^{-1}(0)\right\)+
%c_1(\gamma)\sigma^2 \nonumber \\
c_0(\gamma)\left(||X\bbeta||^2+\sigma^2(1+1/\gamma)L_0
\right)+c_1(\gamma)\sigma^2 \nonumber \\
& \leq & c_2(\gamma)L^*(||X\bbeta||^2+\sigma^2) \label{eq:L0}
\end{eqnarray}
Combining (\ref{eq:Lk1}) and (\ref{eq:L0}) completes the proof.
\newline $\Box$

\subsection{Proof of Theorem \ref{th:upper}}
For all $p_0 \leq r$, Theorem \ref{th:bound} and 
(\ref{eq:pen1}) under the assumption $\pi(r)\geq e^{-c r}$ imply
\begin{eqnarray}
\sup_{\bbeta: ||\bbeta||_0 \leq p_0}
E||X\hat{\bbeta}_{\hat M}-X\bbeta||^2 & \leq & 
\sup_{\bbeta: ||\bbeta||_0 \leq r}
E||X\hat{\bbeta}_{\hat M}-X\bbeta||^2\; 
\leq Pen(r)+c_1(\gamma)\sigma^2 \nonumber \\
& \leq & C_1(\gamma)\sigma^2 r \label{eq:t1.0}
\end{eqnarray}

On the other hand,
applying the general upper bound for the risk of MAP model selector established in
Theorem \ref{th:bound} for models of size $p_0<r$ we have
\begin{equation}
\sup_{\bbeta: ||\bbeta||_0 \leq p_0}
||X\hat{\bbeta}_{\hat M}-X\bbeta||^2 \leq c_0(\gamma)
2\sigma^2(1+1/\gamma)\left(\ln\left\{{p \choose p_0}\pi^{-1}(p_0)\right\}+
\frac{p_0}{2}\ln(1+\gamma)\right)+
c_1(\gamma)\sigma^2 \label{eq:t1.1}
\end{equation}
Abramovich {\em et al.} (2010, Lemma 1) showed that 
${p \choose p_0} \leq (pe/p_0)^{p_0}$. Hence, under the conditions on 
$\pi(p_0)$ in Theorem \ref{th:upper}, for $p_0=1,...,r-1$, (\ref{eq:t1.1}) yields
\begin{eqnarray*}
\sup_{\bbeta: ||\bbeta||_0 \leq p_0}
E||X\hat{\bbeta}_{\hat M}-X\bbeta||^2 & \leq &
c_0(\gamma)2\sigma^2(1+1/\gamma)\left\{(c+1)p_0 \ln(pe/p_0)+\frac{p_0}{2}\ln(1+\gamma) \right\}+c_1(\gamma)\sigma^2 \\
& \leq & C_1(\gamma)\sigma^2 p_0(\ln(p/p_0)+1)
\end{eqnarray*}
Finally, note that for $p_0=r$, as we have already established in (\ref{eq:t1.0}), 
$$
\sup_{\bbeta: ||\bbeta||_0 \leq r}
E||X\hat{\bbeta}_{\hat M}-X\bbeta||^2\;\leq\; C_1(\gamma)\sigma^2 r\;\leq\;
C_1(\gamma)\sigma^2 r(\ln(p/r)+1)
$$
\newline $\Box$

\subsection{Proof of Theorem \ref{th:lower}}
%The proof for the lower bound is based on the approach similar to that in
%the analogous proofs in Birg\'e \& Massart (2001, Theorem 5) and Bunea, 
%Tsybakov \& Wegkamp
%(2007, Theorem 5.1) but modified for a general non-orthogonal design
%and possibility of $p>n$. 
The core of the proof is to find a subset ${\cal B}_{p_0}$ of vectors $\bbeta$, 
where $||\bbeta||_0 \leq p_0$, and the
corresponding subset of mean vectors 
${\cal G}_{p_0}=\{{\bf g} \in \mathbb{R}^n: {\bf g}=X\bbeta,\;\bbeta \in 
{\cal B}_{p_0}\}$ such that for any ${\bf g}_1,\;{\bf g}_2 \in {\cal G}_{p_0}$, 
$||{\bf g}_1-{\bf g}_2||^2 \geq 4s^2(p_0)$ and the Kullback-Leibler divergence
$K(\mathbb{P}_{\bg_1},\mathbb{P}_{\bg_2})=\frac{||\bg_1-\bg_2||^2}{2\sigma^2} \leq (1/16)
\ln {\rm card}({\cal G}_{p_0})$. 
Lemma A.1 of Bunea, Tsybakov \& Wegkamp (2007) will imply then that
$s^2(p_0)$ is the minimax lower bound over ${\cal M}_{p_0}$.

To construct the desired subsets ${\cal B}_{p_0}$ and ${\cal G}_{p_0}$ we consider three
possible cases.

%where evidently ${\rm card}({\cal G}_{p_0})={\rm card}(\tilde{\cal B}_{p_0})$. 

\medskip
\noindent {\sc Case 1}. $p_0 \leq r/2$ 
\newline
Define the subset $\tilde{\cal B}_{p_0}$ of all vectors $\bbeta \in \mathbb{R}^p$
that have $p_0$ entries equal to $C_{p_0}$ defined later,
while the remaining  entries are zeros:
$\tilde{\cal B}_{p_0}=\{\bbeta: \bbeta \in \{\{0,C_{p_0}\}^{p}\},\;
||\bbeta||_0 = p_0 \}$. For $p_0 \leq r/2$
from Lemma 8.3 of Rigollet \& Tsybakov (2010),
there exists a subset
${\cal B}_{p_0} \subset \tilde{\cal B}_{p_0}$ such that for some constant
$\tilde{c}>0$, $\ln{\rm card}({\cal B}_{p_0}) \geq 
\tilde{c} p_0(\ln(p/p_0)+1)$, and for any pair $\bbeta_1,\;\bbeta_2 
\in {\cal B}_{p_0}$, the Hamming distance 
$\rho(\bbeta_1,\bbeta_2)=\sum_{j=1}^p \mathbb{I}\{\bbeta_{1j} \neq \bbeta_{2j}\}
\geq \tilde{c}p_0$.

Consider the corresponding subset of mean vectors ${\cal G}_{p_0}$, 
where ${\rm card}({\cal G}_{p_0})={\rm card}({\cal B}_{p_0})$. 
%${\cal G}_{p_0}=\{{\bf g} \in \mathbb{R}^n: {\bf g}=X\bbeta,\;\bbeta \in 
%\tilde{\cal B}_{p_0}\}$, where evidently ${\rm card}({\cal G}_{p_0})={\rm card}(\tilde{\cal B}_{p_0})$. 
For any $\bg_1,\; \bg_2 \in {\cal G}_{p_0}$ and the
corresponding $\bbeta_1,\;\bbeta_2 \in {\cal B}_{p_0}$ we then have
\begin{equation}
||\bg_1-\bg_2||^2 = ||X(\bbeta_1-\bbeta_2)||^2 \geq \phi_{min}[2p_0] \;
||\bbeta_1-\bbeta_2||^2 \geq \tilde{c}\phi_{min}[2p_0]C^2_{p_0}\;p_0
\label{eq:t4.1}
\end{equation}

On the other hand, by similar arguments, the Kullback-Leibler divergence 
satisfies 
\begin{equation}
K(\mathbb{P}_{\bg_1},\mathbb{P}_{\bg_2})
%\leq \frac{\phi_{max}[2p_0 \wedge r]||\bbeta_1-\bbeta_2||^2}{2\sigma^2}
\leq
\frac{\phi_{max}[2p_0] C^2_{p_0} \rho(\bbeta_1,\bbeta_2)}{2\sigma^2}
\leq \frac{\phi_{max}[2p_0] C^2_{p_0} p_0}{\sigma^2}
\label{eq:t4.2}
\end{equation}
Set now $C^2_{p_0}=(1/16)\sigma^2\tilde{c}(\ln(p/p_0)+1)/\phi_{max}[2p_0]$
and $s^2(p_0)=(1/64)\sigma^2 \tilde{c}^2\tau[2p_0]p_0(\ln(p/p_0)+1)$. 
Then, (\ref{eq:t4.1}) and (\ref{eq:t4.2}) yield
$||\bg_1-\bg_2||^2 \geq 4s^2(p_0)$,  
$K(\mathbb{P}_{\bg_1},\mathbb{P}_{\bg_2}) \leq 
(1/16)\ln {\rm card}({\cal G}_{p_0})$, and Lemma A.1 of Bunea, Tsybakov \& 
Wegkamp (2007) completes the proof.

\medskip
\noindent {\sc Case 2}. $r/2 \leq p_0 \leq r,\;p_0 \geq 8$
\newline
In this case consider the subset
$\tilde{\cal B}_{p_0}=\{\bbeta \in \mathbb{R}^p: \bbeta \in \{\{0,C_{p_0}\}^{p_0},0,...,0\}$ 
and apply Varshamov-Gilbert bound 
(see, e.g. Tsybakov, 2009, Lemma 2.9).
It guarantees the existence of a subset
${\cal B}_{p_0} \subset \tilde{\cal B}_{p_0}$ such that 
$\ln{\rm card}({\cal B}_{p_0}) \geq (p_0/8) \ln 2$ and the Hamming distance
$\rho(\bbeta_1,\bbeta_2) \geq p_0/8$ for any pair $\bbeta_1,\;\bbeta_2 
\in {\cal B}_{p_0}$.  

Note also that for any $\beta_1,\;\bbeta_2 \in {\cal B}_{p_0}$,
$\bbeta_1-\bbeta_2$ has at most $p_0$ non-zero componens and
repeating the arguments for the Case 1, one achieves the minimax lower bound
$s^2(p_0)=C\sigma^2\tau[p_0]p_0 \geq (C/2)\sigma^2\tau[p_0]r$. 

\medskip
\noindent {\sc Case 3}. $r/2 \leq p_0 \leq r,\;2 \leq p_0 < 8$
\newline
For this case, obviously, $2 \leq r < 16$. 
Consider a trivial subset ${\cal B}_{p_0}$ containing
just two vectors $\bbeta_1 \equiv 0$ and $\bbeta_2$ that has $p_0$ nonzero
entries equal to $C^2_{p_0}= (1/64)\sigma^2 \ln 2/\phi_{max}[p_0]$.
For the corresponding mean vectors $\bg_1=X\bbeta_1 = {\bf 0}$ and $\bg_2=X\bbeta_2$,
following (\ref{eq:t4.1}) and (\ref{eq:t4.2}) one has
$$
K(\mathbb{P}_{\bg_1},\mathbb{P}_{\bg_2}) 
\leq \frac{\phi_{max}[p_0] 8 C_{p_0}^2}{2\sigma^2} = 
(1/16) \ln{\rm card}({\cal G}_{p_0})
$$
and 
$$
||\bg_1-\bg_2||^2 \geq \phi_{min}[p_0] p_0 C^2_{p_0} = C\sigma^2\tau[p_0]p_0 \geq 
(C/2) \sigma^2\tau[p_0]r
$$
Applying Lemma A.1 of Bunea, Tsybakov \& Wegkamp (2007) completes the proof.
\newline $\Box$

\subsection{Proof of Theorem \ref{th:nonorth}}
%Let $M_{p_0} \in {\cal M}_{p_0}$ be the true (unknown) model in (\ref{eq:model1}).
We want to show that under the conditions of Theorem \ref{th:nonorth} we can
reduce the rate in the upper bound (\ref{eq:upper}) for the risk of the MAP 
model selector
established in Theorem \ref{th:upper} over ${\cal M}_{p_0}$ by the factor 
$\tau_p[2p_0]$. Recall that we derived (\ref{eq:upper}) from the general upper bound
(\ref{eq:bound}) in Theorem \ref{th:bound} by considering models $M$ of size $p_0$. 
Consider now $1 \leq \kappa_{p1} \leq p_0 \leq \kappa_{p2} \leq r/2$ and apply 
(\ref{eq:bound}) for models $M'$ of less size 
$p'_0=\lceil \tau_p[2p_0]p_0 \rceil \geq 1$. 
Consider an arbitrary $\bbeta_p$ with $||\bbeta_p||_0 = p_0$.
Under the conditions of Theorem \ref{th:nonorth} on the prior $\pi_p(\cdot)$,
(\ref{eq:bound}) implies 
\begin{equation}
E||X_p\hat{\bbeta}_{{\hat M}_p}-X_p\bbeta_p||^2
\leq C_1(\gamma)\inf_{M'}
||X_p\bbeta_p-X_p\bbeta_{pM'}||^2
+C_2(\gamma)\sigma^2p'_0(\ln(p/p'_0)+1),  \label{eq:t5.0}
\end{equation}
where $X_p\bbeta_{pM'}$ is the projection of the mean vector $X_p\bbeta_p$
on the span of $M'$.
Comparing (\ref{eq:upper}) and (\ref{eq:t5.0}) illustrates that
reduction of a model size introduces the bias.  
On the other hand, under Assumptions (D.1) and (D.2), a straightforward calculus
shows then that the variance term decreases to the desired order 
$p_0 \tau_p[2p_0](\ln(p/p_0)+1)$. 
The idea of the proof will be based on finding  a model $M'_*$
such that the resulting bias term will be at most of the same order as the reduced variance. 

Consider the model $M$ of size $p_0$ corresponding to $\bbeta_p$
and any of its submodels $M'$ of size $p'_0$ defined above. Then, $M/M'$ is
evidently a subset of predictors from $M$ not included in $M'$ and
$D_{M/M'}=D_M-D_{M'}$, where
diagonal indicator matrices $D$'s were introduced in Section \ref{sec:MAP}.
%$D_{M/M',jj}=\mathbb{I}\{\beta_j \in M, \beta_j \notin M'\}$.
By straightforward calculus one then has
\begin{eqnarray}
||X_p\bbeta_p-X_p\bbeta_{pM'}||^2 & = & 
\bbeta'_pD_{M/M'}(D_{M/M'}(D_MX_p'X_pD_M)^{+}D_{M/M'})^{+}D_{M/M'}\bbeta_p  \nonumber \\
& = & \bbeta'_p G_{p,{\cal J}_{M/M'}} \Lambda^{-1}_p(M,M')G'_{p,{\cal J}_{M/M'}} \bbeta_p 
\;\leq \; \tilde{\phi}_{M,M'}^{-1} ||G'_{p,{\cal J}_{M/M'}}\bbeta_p||^2 \nonumber \\
& \leq & \tilde{\phi}_{M,M'}^{-1} ||\bbeta_p||^2_{\infty} (p_0-p'_0), \label{eq:t5.1}
\end{eqnarray}  
where the matrices $G$ and $\Lambda$ and the
minimal eigenvalue $\tilde{\phi}_{M,M'}$
were defined in Section \ref{subsec:multi}.

Among all submodels $M' \subset M$ of size $p'_0$, choose $M'_*$ with the 
maximal $\tilde{\phi}_{M,M'}$. 
Then, $\tilde{\phi}_{M,M'_*} \geq \tilde{\phi}_p[p_0]$ and 
(\ref{eq:t5.1}) and Assumption (B) yield 
$$
||X_p\bbeta_p-X_p\bbeta_{pM'_*}||^2 \leq C \tau_p[2p_0]p_0\ln((p/p_0)+1)
$$
Hence, we proved that under assumptions on the prior, Assumptions (D.1, D.2) 
and (B), the upper bound for the
risk of the MAP model selector over ${\cal M}_{p_0}$ is of the minimax order 
$\tau_p[2p_0]p_0\ln(p/p_0)+1)$. Assumption (D.3) guarantees
that the ``least-favorable'' sets 
${\cal B}_{p_0}$ constructed in the proof of Theorem \ref{th:lower} 
satisfy the additional Assumption (B) on $\bbeta_p$ and, therefore,
the minimax lower bound (\ref{eq:lower}) is not reduced.
\newline $\Box$

\end{document}